\documentclass[12pt]{article}

\usepackage{amsmath,amsthm}
\usepackage{amsfonts,xcolor,graphicx,url}
\usepackage{epstopdf,epsfig}

\usepackage[margin=2.5cm]{geometry}

\usepackage{titlesec}
\titleformat{\section}
{\normalfont\bfseries}{\thesection}{1em}{}
\titleformat{\subsection}
{\normalfont}{\thesubsection}{1em}{}

\numberwithin{equation}{section}

\theoremstyle{plain}

\theoremstyle{definition}

\theoremstyle{remark}
\newtheorem{remark}{Remark}

\begin{document}

\title{
Piecewise flat Ricci flow of compact without boundary three-manifolds
}

\author{Rory Conboye}

\date{}
\maketitle

\let\thefootnote\relax\footnotetext{
	\hspace{-0.75cm}
	Department of Mathematics and Statistics
	\hfill Present Address: 
	\\
	American University
	\hfill Department of Mathematics 
	\\
	4400 Massachusetts Avenue, NW
	\hfill Munster Technological University
	\\
	Washington, DC 20016, USA
	\hfill Bishopstown, Cork T12 P928, Ireland
	\\
	${}^{} $
	\hfill Rory.Conboye@mtu.ie
}
\let\thefootnote\svthefootnote

\begin{abstract}

Using a recently developed piecewise flat method, numerical evolutions of the Ricci flow are computed for a number of manifolds, using a number of different mesh types, and shown to converge to the expected smooth behaviour as the mesh resolution is increased.
The manifolds were chosen to have varying degrees of homogeneity, and include Nil and Gowdy manifolds, a three-torus initially embedded in Euclidean four-space, and a perturbation of a flat three-torus.
The piecewise flat Ricci flow of the first two are shown to converge to known smooth Ricci flow solutions, with the remaining two flowing asymptotically to flat metrics.

\


\


\end{abstract}

\section{Introduction}


The Ricci flow is a uniformizing flow on manifolds, acting like a heat flow. It was developed by Richard Hamilton in the 1980's \cite{HamRF} as an approach to proving the Thurston geometrization conjecture, and therefore the Pioncar\'{e} conjecture, an endeavour that was eventually successful with the work of Grisha Pereleman in the early 2000's \cite{Pereleman1,Pereleman2}. Since then, Ricci flow has remained an important tool for investigating the interplay between geometry and topology, and has also found a number of applications in image analysis, from facial recognition \cite{GuZengFacial} to cancer detection \cite{GuCancer}, and in space-time physics \cite{WoolgarRFphys,WiseRF-BH,WiseRF-CFT}. Applications such as these require a robust numerical approach, and preliminary work has even shown the potential for numerical simulations to inform analytical Ricci flow research \cite{GarfIsen}.


Two dimensional numerical Ricci flow has seen a lot of progress recently, from Chow and Luo's combinatorial Ricci flow \cite{CombRF} to the discrete surface Ricci flow of Gu et al. \cite{GuDiscSurfRF,GuSurveyDRF}. These approaches use piecewise linear approximations of smooth surfaces with the Ricci flow given by a conformal deformation, utilizing the equivalence of the Ricci and scalar curvatures in two dimensions. In three dimensions, Garfinkle \& Isenberg \cite{GarfIsen,GarfIsen08} have used finite difference methods to find the critical behaviour for a one-parameter family of spherically symmetric neck pinch geometries on $S^3$. Crucially, this was done before the problem was completely understood analytically. A three-dimensional piecewise flat approach was also introduced a number of years ago \cite{SRF}, based on intuitions gained from Regge calculus \cite{Regge}, and applied to the same neck pinch geometries \cite{MNeckPinch1}. However, its effectiveness seems to be restricted to highly symmetric geometries, with its computational mesh adapted to these symmetries.


The computations in this paper are based on a new piecewise flat approach developed in \cite{PLCurv}. Significantly, this approach has no symmetry restrictions, a particularly broad freedom in the choice of mesh, and gives results that are independent of such mesh choices to reasonable levels of precision. 
The Ricci flow evolution consists of a set of independent equations for each edge in the graph, giving it the potential to be highly parallelizable. 
The equations depend on a piecewise flat Ricci curvature, given in a computationally efficient combinatorial form, which has already been shown to converge to the smooth Ricci curvature with increasing mesh resolutions in \cite{PLCurv}, for a number of different manifolds and mesh types.
There are also substantial advantages in the use of piecewise flat manifolds, including the retaining of a manifold structure after discretization, avoiding issues with coordinate singularities due to the coordinate-independent structure, and having the topology completely determined by the piecewise flat graph.


Evolutions of this piecewise flat Ricci flow are computed for a number of different situations, including:
\begin{enumerate}
	\item four different manifolds with different levels of homogeneity;
	
	\item three different types of triangulation mesh, each with different orientations of edges; 
	
	\item three different mesh resolutions for each triangulation type.
\end{enumerate}
The results of these computations show that there is no dependence on particular manifold symmetries, or on specially adapted triangulations. The computations also show that equivalent results are given for different triangulation types, up to a reasonable level of precision, and that the behaviour for higher mesh resolutions converges to the expected smooth behaviour. 
To ensure that the effects of the different triangulation types and mesh resolutions can be compared directly, the equations are evolved using a simple Euler method with a consistent number and size of time steps for each manifold.


The four manifolds were chosen to have a variety of different characteristics, but all are compact, with topologies close to a three-torus. The first manifold is a Nil geometry, one of Thurston's eight homogeneous geometries, with analytic solutions for the normalized and non-normalized Ricci flow found by Isenberg \& Jackson \cite{IsenJackRFHomog} and Knopf \& McLeod \cite{KnMcModel} respectively.
The second manifold is a Gowdy three-geometry, one of the first non-positive-definite curvature geometries shown to have a convergent Ricci flow by Carfora, Isenberg \& Jackson \cite{CarfIsenGowdy}. This has a two-dimensional isometry group, reducing the Ricci flow equations to a set of two coupled partial differential equations (PDEs) for which a numerical solution can easily be found. 
A three-torus initially embedded in Euclidean four-space was chosen for the third manifold, having a one-dimensional isometry group. For the last, a perturbation of a flat three-torus with no continuous isometries was chosen. The latter two manifolds are expected to Ricci flow asymptotically to a flat three-torus.


The rest of the paper begins with an introduction to the new piecewise flat Ricci flow method, summarising the main results and details from \cite{PLCurv}. This is followed by a description of the three different types of triangulation. Each of the four manifolds is then dealt with separately in sections \ref{sec:Nil} to \ref{sec:Pert}, with subsections to introduce the smooth manifolds, provide specific details about the triangulations, and both display and discuss the results of the Ricci flow computations.

\section{Ricci flow and triangulations}

\subsection{Piecewise flat Ricci flow}
\label{sec:PFRF}


Simplicial piecewise flat manifolds in three dimensions are formed by joining Euclidean tetrahedra together, identifying the triangular faces of neighbouring segments. The topology of such a manifold is completely determined by the resulting graph, and the geometry by the set of edge-lengths. 
While each pair of neighbouring tetrahedra form a consistent Euclidean space, the dihedral angles of the tetrahedra meeting at a given edge may not necessarily sum to $2 \pi$ radians, with the difference known as the deficit angle
\begin{equation}
\epsilon_\ell := 2 \pi - \sum_t \theta_t,
\end{equation}
for the dihedral angles $\theta_t$ at the edge $\ell$. A smooth manifold can then be approximated by first constructing a tetrahedral lattice, using geodesic segments as edges, and then defining a piecewise flat manifold using the same graph, with the edge-lengths defined by the lengths of the corresponding geodesic segments. A good approximation will have uniformly small deficit angles, which can be achieved by using a lattice with high resolution in areas of high curvature.

\begin{figure}[h]
	\centering
	\includegraphics[scale=1]{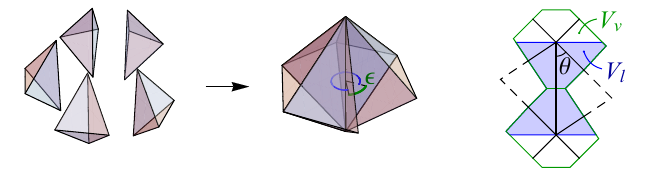}
	\caption{The deficit angle at an edge, and a cross section of the region around an edge showing the vertex and edge volumes.}
	\label{fig:tet}
\end{figure}


On a smooth manifold $M$, the Ricci flow changes the metric $g$ according to the Ricci curvature $Rc$ at each point, with an optional second term giving a normalized Ricci flow,
\begin{equation}
\frac{d g}{d t} = - 2 Rc \, ,
\qquad {\text or} \qquad
\frac{d g}{d t} = - 2 Rc + \frac{2}{3} \widetilde{R}_M \, g \, ,
\end{equation}
where the volume can be kept constant with the use of the average scalar curvature $\widetilde{R}$ of the manifold. 
Both equations tend to reduce the value of the metric components associated with high positive Ricci curvature and increase those associated with high negative Ricci curvature, with the rate depending on the strength of the curvature. It was shown in \cite{PLCurv} that the effect of the Ricci flow on the lengths of geodesic segments can be given by the integral of the Ricci curvature along and  tangent to the segment. This naturally leads to a piecewise flat approximation of the smooth Ricci flow as a set of independent equations for the fractional change in the edge-lengths,
\begin{equation}
\label{eq:PFRF}
\frac{1}{|\ell|} \frac{d |\ell|}{d t}
 = - Rc_\ell
 \, ,
\qquad
\frac{1}{|\ell|} \frac{d |\ell|}{d t}
= - Rc_\ell + \frac{1}{3} \widetilde{R}_S
\, ,
\end{equation}
for piecewise flat approximations  $Rc_\ell$, of the curvature tangent to the corresponding geodesic segments, and $\widetilde{R}_S$, of the average scalar curvature of the manifold. These equations will continue to give a good approximation for the smooth Ricci flow as long as the lattice gives uniformly small deficit angles.


All that remains is to give the piecewise flat curvature approximations $Rc_\ell$ and $\widetilde{R}_S$. The first is given in terms of the scalar curvature $R_v$ at the vertices $v_1$ and $v_2$ bounding $\ell$, and the sectional curvature $K_\ell^\perp$ orthogonal to the edge $\ell$. 
These are in turn defined over the volumes $V_v$ of a dual tessellation, with barycentric duals used in this paper, and edge-volumes $V_\ell$, defined as the union of the vertex volumes on either end of the edge, bounded by surfaces orthogonal to $\ell$ at each vertex. A cross-section of these volumes can be seen in figure \ref{fig:tet}. From \cite{PLCurv},
\begin{subequations}
\label{eq:RandK}
\begin{align}
Rc_\ell
 &= \frac{1}{4} (R_{v_1} + R_{v_2}) - K_\ell^\perp \, ,
 \label{eq:Rc_l} \\
K_\ell^\perp
 &= \frac{1}{V_\ell} \left(
 |\ell| \, \epsilon_\ell + \sum_i \frac{1}{2} |\ell_i| \cos^2 (\theta_i) \epsilon_i
 \right) ,
 \label{eq:K_l} \\
R_v 
 &= \frac{1}{V_v} \sum_{i} |\ell_i| \, \epsilon_i \, ,
\qquad
\widetilde{R}_S 
  = \frac{2}{V_S} \sum_{i} |\ell_i| \, \epsilon_i \, .
 \label{eq:R_v}
\end{align}
\end{subequations}
The indices $i$ correspond to the edges intersecting the volumes $V_v$ and $V_\ell$ (aside from $\ell$ itself), with $\epsilon_i$ representing the deficit angles and $\theta_i$ the angle between $\ell_i$ and the edge $\ell$. 
The expressions for $K_\ell^\perp$ and $R_v$ were found by constructing volume-integrals of each curvature, with the deficit angles shown to represent local surface-integrals of the sectional curvature, and then divided by the volume to give an average. The expression for $\widetilde{R}_S$ can also be seen to come directly from the Regge action \cite{Regge}, divided by the volume $V_S$ of the piecewise flat manifold.


Computations in \cite{PLCurv} have successfully shown these expressions to converge to their corresponding smooth curvature values as the mesh resolution is increased, for a variety of different manifolds, and provide good approximations for reasonably dense meshes. The computations in this paper should provide even more support for the effectiveness of these piecewise flat curvature expressions. Similar expressions have also been developed for the extrinsic curvature on piecewise flat manifolds using a similar approach \cite{PLExCurv}, with computations also showing good approximation of, and convergence to, their corresponding smooth curvatures.

\subsection{Triangulation types}
\label{sec:Tri}

\begin{figure}[h]
	\centering
	\includegraphics[scale=1.2]{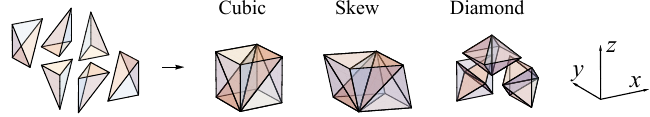}
	\caption{The three different block types, with the six tetrahedra of the cubic block on the far left, and a slight separation of the three diamond shapes forming the diamond block.}
	\label{fig:blocks}
\end{figure}

For the purposes of this paper, a triangulation-type is defined as a family of simplicial graphs on a manifold, where the number of tetrahedra can be scaled in some regular way by choosing different members of the family. This scaling provides control over mesh refinements and is used to show the convergence of piecewise flat constructions to their smooth counterparts. In this section, three such triangulation types will be defined on a three-torus ($T^3$) topology, with a fundamental domain of size $[0,1]^3$ for a set of coordinates $x$, $y$, $z$. The three different triangulation types provide a variety of edge orientations, ensuring that computation results are not based on favourable alignments alone, and provide differing levels of suitability for different dual tessellations. The building blocks for each triangulation type are defined below, and shown in figure \ref{fig:blocks}.

\begin{enumerate}
	\item The \emph{cubic} block is cube-shaped and composed of six tetrahedra. There are seven edges for a given vertex, three along the $x$, $y$ and $z$ coordinates, three face-diagonals and a body-diagonal, with the tetrahedra specified so that the orientation of the face-diagonals on opposite sides agree. This is the most simple construction but it only borders on being a Delaunay triangulation for a flat metric, with the circumcenters of all six tetrahedra coinciding at the centre of the cube.
	
	\item A \emph{skew} block has the same structure as the cubic block but with the vertices $v_x$ and $v_z$ skewed so that $v_x = (1, -1/3,0)$ and $v_z = (-1/3,-2/9,1)$. This block forms a different tiling of the same domain, and gives a strongly Delaunay triangulation for a flat metric.
	
	\item A \emph{diamond} block is constructed by forming a set of four tetrahedra in a diamond shape around each coordinate edge, with the edges in the outer ring parallel with the other coordinate directions. This block contains two distinct vertices, instead of the single vertex for the other two blocks, and 14 edges, six that are parallel with the three coordinates and eight forming a set of four body-diagonal-type edges.
\end{enumerate}

In each case, a single block can be used as a simplicial graph of $T^3$ by identifying faces on opposite sides. A regular grid of blocks can then be used to refine the mesh, with faces on opposite sides of the grid identified. Once the graph has been defined the vertices can be labelled, with each edge, triangle and tetrahedron then specified by the list of vertices in its closure. Ambiguities in labelling can arise for grids that have less than three vertices in a single direction, for example, all vertices in a single cubic block will have the same label. In these cases the grid is duplicated to form a slightly larger covering space before labelling. The computations are still performed for the original domain alone, with any geometric data copied to the covering space.

Unfortunately, the piecewise flat Ricci flow equations are not stable when applied directly to the cubic and skew triangulation-types, seemingly due to an overdetermination of the system. However, this instability can be avoided by treating the interior of each block as flat, see \cite{PLRFinstab}. In practice this is done by re-defining the body-diagonal edge-lengths so that the deficit angles at these edges are essentially zero. The piecewise flat Ricci flow is still performed on a simplicial triangulation, where dihedral angles and volumes can be computed in a standard way, but with the lengths of the body-diagonals adjusted for each time-step to give an essentially flat interior for each block. This approach has led to stable evolutions of the piecewise flat Ricci flow for all of the computations in this article. Details of the instability and proofs of the stabilizing method are provided in a related paper \cite{PLRFinstab}.

\

\begin{remark}
The true advantage of piecewise flat approximations is the complete encoding of the topology in the simplicial graph. In particular, for compact without boundary manifolds, once the boundary identifications have been made and the labels assigned, the boundaries vanish. There are no boundary conditions required for functions on these piecewise flat manifolds, since the graph already has the appropriate topology. There is also no need for any coordinate systems, once geodesic lengths have been computed they just become properties of the graph. As well as avoiding coordinate singularities or translations between multiple charts, the geometry is also clear from the lengths of the edges and does not need to be separated from the characteristics of a coordinate system. This even led Tullio Regge to title his 1961 paper \cite{Regge} as ``General relativity without coordinates''.
\end{remark}

\section{Nil manifolds}
\label{sec:Nil}

\subsection{Smooth manifold}


Thurston's geometrization conjecture states that any closed 3-manifold can be decomposed into a set of irreducible parts, each of which admits a locally homogeneous geometry. The universal covers of these geometries form a set of eight model geometries, one of which is the Nil geometry, the geometry of the continuous Heisenberg group. The Heisenberg group is a nilpotent Lie group, and can be represented by the group of $3 \times 3$ matrices of the form
\begin{equation}
\label{eq:Hgp}
\left(
\begin{array}{ccc}
1 & x & z \\
0 & 1 & y \\
0 & 0 & 1 
\end{array}
\right)
\end{equation}
with real number entries $x$, $y$ and $z$, which can then be seen as a set of coordinates on $\mathbb{R}^3$.


Analytic solutions for the normalized Ricci flow of the Nil geometry were first given by Isenberg and Jackson \cite{IsenJackRFHomog}, with the quasi-convergence of the non-normalized Ricci flow later studied by Knopf and McLeod \cite{KnMcModel}. These approaches make use of work by Milnor \cite{Milnor76} and by Ryan and Shepley \cite{RyanShepley75}, where an orthogonal frame of 1-forms $\{ \theta^a \}$ can be found in which the metric is diagonal
\begin{equation}
\label{eq:nilDiag}
g = A (\theta^1)^2 + B (\theta^2)^2 + C (\theta^3)^2 ,
\end{equation}
and the structure constants $C^a_{b d}$ (such that $d \theta^a = C^a_{b d} \theta^b \wedge \theta^d$) are all zero except for $C^1_{2 3} = - C^1_{3 2} = \lambda$, for some constant $\lambda$. The Ricci curvature is also diagonal in this frame, which reduces both the normalized and non-normalized Ricci flows to a set of ordinary differential equations for the metric components $A$, $B$ and $C$, with the frame and structure constants invariant to the flows. The constant $\lambda$ can be seen as a scaling factor on the $z$ coordinate in the matrix representation of the Heisenberg group (\ref{eq:Hgp}),
\begin{equation}
\label{eq:Hgp2}
\left(
\begin{array}{ccc}
1 & x & z/\lambda \\
0 & 1 & y \\
0 & 0 & 1 
\end{array}
\right).
\end{equation}
In \cite{IsenJackRFHomog} this constant is taken as the volume element ($\lambda = \sqrt{A B C}$), coming from a weighted Levi-Civita antisymmetric tensor in equation (4), which is invariant under the volume-preserving normalized Ricci flow. Solutions for $A$, $B$ and $C$ are then given in equations (29). A value of $-2$ is used for $\lambda$ in \cite{KnMcModel}, with solutions for the non-normalized Ricci flow given in equations (15). The latter solutions can also be found in the Ricci flow textbooks \cite{ChowKnopfRF,ChowLuNiRF}, in equations (1.8) and (4.62) respectively.

It can easily be shown that the scalar curvature $R = -\frac{\lambda^2 A}{2 B C}$, with solutions for the metric functions under the normalized Ricci flow then given by the equations:
\begin{equation}
\label{eq:NilSolnIJ}
A = A_0 \left(1 - \frac{16}{3} R_0 t\right)^{-1/2}, \quad
B = B_0 \left(1 - \frac{16}{3} R_0 t\right)^{1/4}, \quad
C = C_0 \left(1 - \frac{16}{3} R_0 t\right)^{1/4},
\end{equation}
and under the non-normalized Ricci flow by:
\begin{equation}
\label{eq:NilSolnKM}
A = A_0 (1 - 6 R_0 t)^{-1/3}, \quad
B = B_0 (1 - 6 R_0 t)^{1/3}, \quad
C = C_0 (1 - 6 R_0 t)^{1/3},
\end{equation}
for initial values $A_0$, $B_0$ and $C_0$, and the initial scalar curvature $R_0 = -\frac{\lambda^2 A_0}{2 B_0 C_0}$. The normalized Ricci flow solutions above agree with those of \cite{IsenJackRFHomog} for $\lambda = \sqrt{A B C}$, and the non-normalized solutions with \cite{KnMcModel} for $\lambda = -2$.


For computational purposes, the initial metric (\ref{eq:nilDiag}) is chosen so that $A_0 = B_0 = C_0 = 1$. In the coordinates given by the matrix representation of the Heisenberg group in equation (\ref{eq:Hgp2}), the metric takes the form
\begin{equation}
\label{eq:gNil}
d s^2 = d x^2 + d y^2 + (d z + \lambda \, x \, d y)^2.
\end{equation}
The orthogonal frame is obtained from these coordinates by the identification $\{\theta^1, \theta^2, \theta^3\} = \{d z, dy + \lambda \, x \, d z, d x\}$. In order to obtain a compact computational domain, the continuous Heisenberg group can be quotiented by the finite Heisenberg group, with integer entries in the upper-right triangle. This gives a fundamental domain with the $x$ and $y$ coordinates having a range of $[0, 1]$ and the $z$ coordinate a range of $[0, \lambda]$. The $yz$ planes have a two-torus topology, while the identification of the $yz$ faces at $x=0$ and $x=1$ requires a sort of \emph{twist}, as shown in figure \ref{fig:NilTri}.

\subsection{Triangulating manifold}

\begin{figure}[h!]
	\begin{center}
		\includegraphics[scale=1.2]{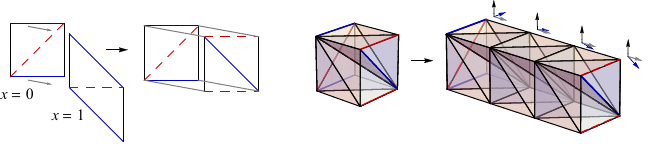}
	\end{center}
	\caption{Relations between the $yz$-faces at $x=0$ and $x=1$ for $\lambda = 1$, and single and three-block triangulations with these face relations. The orthogonal frames for $\{\theta^i\}$ are also indicated for different values of $x$ along the three-block triangulation.}
	\label{fig:NilTri}
\end{figure}

For a single cubic block triangulation, the decomposition into six tetrahedra is done in a slightly different way to section \ref{sec:Tri}, see figure \ref{fig:NilTri}, so that the face-diagonals match the face relations for the fundamental domain. The geometry is isometric in the $y z$-planes, with a standard two-torus topology, so triangulations only require a single block in the $y$ and $z$ directions. Although the Nil geometry is completely homogeneous, the metric representation in (\ref{eq:gNil}) depends on the $x$-coordinate, so increased resolutions are given by joining copies of the standard cubic block along the $x$-direction, with the complete collection of blocks still spanning $[0, 1]$. This has the effect of sub-dividing the \emph{twist} into multiple parts, and reduces the size of the deficit angles. Vertices, edges and triangles on opposite faces are then identified, paying attention to the $yz$-face relations. Adapting the skew and diamond triangulation types to the Nil face relations is not so straight forward, so only the cubic triangulation type will be used here. However, since the $y$-edges are not orthogonal to the $z$-edges beyond $x = 0$, each block in the grid can be considered as a slightly different triangulation type.

In order to compare with the analytic results of Isenberg \& Jackson \cite{IsenJackRFHomog} and Knopf \& McLeod \cite{KnMcModel}, the normalized Ricci flow is computed for triangulations with $\lambda = 1$ and the non-normalized Ricci flow for triangulations with $\lambda = -2$. Grids of one, two and three blocks are then used for both situations. For the $\lambda = 1$ case, these triangulations have a range of $[0,1]$, $[0,1/2]$ and $[0,1/3]$ respectively in the $y$ and $z$ directions, in order to keep the tetrahedra close to regular. The original domain can then be recovered with four and nine copies for the two and three block triangulations. For the $\lambda=-2$ case, the twist is in the opposite direction, so the triangulations are a reflection in the $x y$-plane of those shown in figure \ref{fig:NilTri}. A range of $[0,1/2]$ is used in the $x$-direction here, to keep from stretching the cubes too much, with the $y$ and $z$ directions given ranges of $[0,1/2]$, $[0,1/4]$ and $[0,1/6]$ for the three different grid sizes. Again, the original domain can be recovered with grids of copies with dimensions $2 \times 2 \times 4$, $2 \times 4 \times 8$ and $2 \times 6 \times 12$ in the $x$, $y$ and $z$ directions.

\subsection{Results of evolutions}

\begin{figure}[h!]
	\begin{center}
		\includegraphics[scale=1]{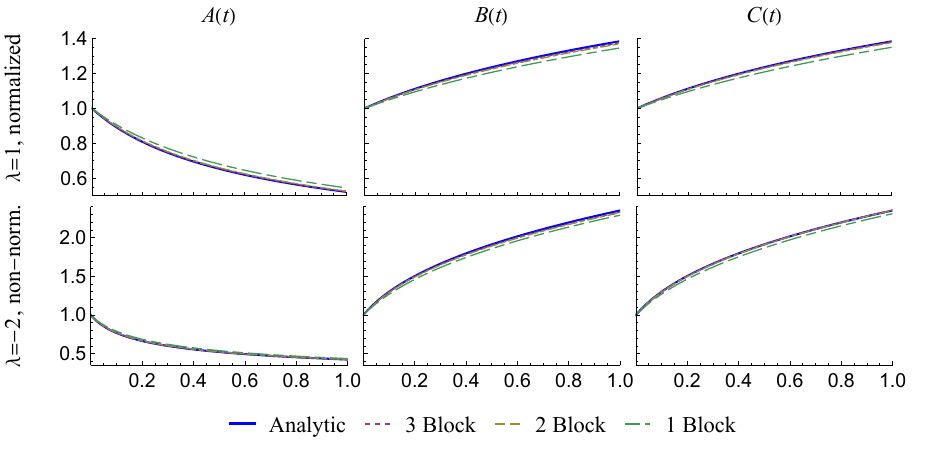}
	\end{center}
	\caption{Graphs of the metric functions in time for both the normalized and non-normalized Ricci flows, all showing convergence to the analytic solutions as the resolution is increased.}
	\label{fig:Nilt}
\end{figure}

Evolutions of the normalized and non-normalized piecewise flat Ricci flow (\ref{eq:PFRF}) are computed using the Euler method, with 200 steps of size 0.005 for each triangulation. Since Nil manifolds are homogeneous, the values of the metric functions can be given by the square of the geodesic lengths for a unit coordinate length in the $x$ and $z$ directions for $C$ and $A$ respectively, and in the $y$-direction at $x=0$ for $B$. The resulting values are graphed in figure \ref{fig:Nilt} as piecewise linear curves, along with the analytic solutions from equations (\ref{eq:NilSolnIJ}) and (\ref{eq:NilSolnKM}). The piecewise flat curves can be seen to give good approximations for the analytic solutions, even for the single block triangulations with only six tetrahedra, and there is a clear convergence to the analytic solutions as the number of blocks is increased.

The piecewise flat values for $A$, $B$ and $C$ have also been fitted to a general function matching the form of the analytic equations (\ref{eq:NilSolnIJ}) and (\ref{eq:NilSolnKM}),
\begin{equation}
	\label{eq:NilParam}
	f(t) = (1 + a \, t)^b .
\end{equation}
The resulting best-fit parameters for $a$ and $b$ are shown in table \ref{tab:NilParam}, for each of the functions and both evolutions, with $R$-squared values of greater than $0.999999$ in all cases. These again show a clear convergence to the analytic solutions from \cite{IsenJackRFHomog} and \cite{KnMcModel}, as the resolution is increased, with very close approximations for the three-block triangulations in particular.

\begin{table}[h!]
	\centering
	\begin{tabular}{lcc|cc|cc
		}
		&\multicolumn{2}{c|}{$A$}
		&\multicolumn{2}{|c|}{$B$}
		&\multicolumn{2}{|c}{$C$}
		\\
		\hline
		Normalized 
		& a & b
		& a & b
		& a & b
		\\
		1-Block 
		& $2.043$ & $-0.545$ 
		& $2.065$ & $0.265$ 
		& $2.033$ & $0.270$ 
		\\
		2-Block 
		& $2.461$ & $-0.516$ 
		& $2.523$ & $0.249$ 
		& $2.476$ & $0.256$ 
		\\
		3-Block 
		& $2.575$ & $-0.508$ 
		& $2.632$ & $0.247$ 
		& $2.583$ & $0.253$ 
		\\
		Analytic 
		& $2.667$ & $-0.5$ 
		& $2.667$ & $0.25$ 
		& $2.667$ & $0.25$ 
		\\
		&&&&&& \\
		\hline
		Non-norm. 
		& a & b
		& a & b
		& a & b
		\\
		1-Block 
		& $ 9.52$ & $-0.351$ 
		& $ 9.91$ & $0.347$ 
		& $ 9.75$ & $0.353$ 
		\\
		2-Block 
		& $11.42$ & $-0.339$ 
		& $11.60$ & $0.332$ 
		& $11.45$ & $0.339$ 
		\\
		3-Block 
		& $11.90$ & $-0.335$ 
		& $12.05$ & $0.329$ 
		& $11.83$ & $0.336$ 
		\\
		Analytic 
		& $12$ & $-0.333$ 
		& $12$ & $0.333$ 
		& $12$ & $0.333$ 
		\\
	\end{tabular}
	\caption{The best fit parameters for expressions of the form (\ref{eq:NilParam}), showing a clear convergence to the analytic solutions, with $R$-squared values greater than $0.999999$ in all cases.
		}
	\label{tab:NilParam}
\end{table}

\section{Gowdy manifold}
\label{sec:Gowdy}

\subsection{Smooth manifold}

The Gowdy manifolds are a family of 3-geometries that can be foliated into isometric 2-surfaces and are periodic in the direction orthogonal to these surfaces. These manifolds first appeared as the spatial part of a cosmological space-time model introduced by Robert Gowdy \cite{Gowdy}, with dynamics capable of containing gravitational waves. They were then used by Carfora, Isenberg and Jackson \cite{CarfIsenGowdy} to show that the Ricci flow can be convergent for manifolds with non-positive-definite curvature. With coordinates $x$ and $y$ spanning the isometric 2-surfaces, and $\theta$ orthogonal to them, the metric can be written in the general form
\begin{equation}
d s^2 = e^{f + W} d x^2 + e^{f- W} d y^2 + e^{2 a} d \theta^2 \, ,
\end{equation}
with a constant function $f$, and functions $W$ and $a$ depending on the value of $\theta$. The non-normalized Ricci flow of this metric can be represented by a coupled set of partial differential equations (PDEs) for these functions:
\begin{equation}
\partial_t f = 0 ,
\quad
\partial_t a = \frac{1}{2} e^{-2 a} \left(\partial_\theta W\right)^2 ,
\quad
\partial_t W = e^{-2 a} \left(
\partial_\theta^2 W - \partial_\theta a \cdot \partial_\theta W
\right) .
\end{equation}
These equations were easily solved numerically for initial-time functions
\begin{equation}
f_0 = 0 , \qquad a_0 = 0, \qquad W_0 = 0.1 \sin \theta ,
\end{equation}
requiring nothing more than the built-in PDE solver in \emph{Mathematica}. The normalized Ricci flow equations are not so straight forward to solve, particularly in ensuring that $f$ is constant for each value of $t$. As a result, only the non-normalized Ricci flow is investigated here, though the volume is also close to being invariant for the non-normalized flow. It should be noted that there is no extra difficulty in computing the piecewise flat normalized Ricci flow.

\subsection{Triangulating manifold}

The Gowdy manifold is the simply connected covering space of a compact $T^3$ manifold with fundamental domain $[0,s_x] \times [0,s_y] \times [0, 2\pi]$, for any values $s_x$ and $s_y$. This can be used as a compact domain for the piecewise flat computations. Since the manifold is isometric in each $x y$-plane, triangulations only require a single block in the $x$ and $y$ directions. Different resolutions for cubic and skew triangulation types are provided by grids of $6$, $12$ and $24$ blocks in the $\theta$ direction, with $3$, $6$ and $12$ diamond blocks having the same number of vertices and edges. To keep the tetrahedra close to regular, both $s_x$ and $s_y$ are chosen to be $2$, $1$, $1 \over 2$ and $1 \over 4$ for the $3$, $6$, $12$ and $24$ block grids respectively, with $4$, $16$ and $64$ copies of the latter three grids giving the same domain for all triangulations. The piecewise flat curvatures for the cubic and skew triangulations have already been shown in \cite{PLCurv} for the initial manifold.

\begin{figure}[h!!]
	\begin{center}
		\includegraphics[scale=1]{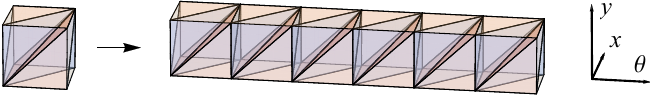}
	\end{center}
	\caption{
	A grid of six cubic blocks along the $\theta$-direction for the lowest resolution triangulation, with the vertices, edges and triangles on opposite sides identified to give a $T^3$ topology.
	}
	\label{fig:GowdyTri}
\end{figure}

\subsection{Results of evolutions}

\begin{figure}[h!]
	\begin{center}
		\includegraphics[scale=1]{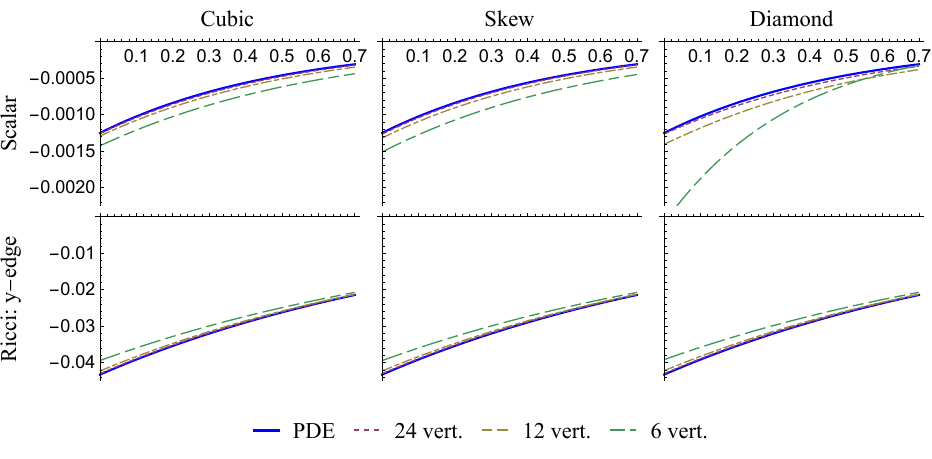}
	\end{center}
	\vspace{-0.5cm}
	\caption{Piecewise-linear graphs of the curvature values in time, at $\theta = \pi/3$, showing a clear convergence of the piecewise flat to the PDE solution as the resolution is increased.}
	\label{fig:Gowdyt}
\end{figure}

The piecewise flat Ricci flow (\ref{eq:PFRF}) was applied to each of the triangulations, using the Euler method with $70$ time-steps of size $0.01$. Piecewise-linear graphs for the flow of the piecewise flat scalar curvature at $\theta=\pi/3$ are shown on the top row of figure \ref{fig:Gowdyt}, together with the scalar curvature from the numerical PDE solution. The bottom row of the figure shows the Ricci curvatures along the $y$-direction at $\theta=\pi/3$. The $\theta$ value of $\pi/3$ was chosen since the initial curvatures are all non-zero here, and there is a vertex at this value of $\theta$ for all of the triangulations. The $y$-edge was selected for the Ricci curvature as this is the only edge that is part of all three triangulation types. The graphs show convergence to the PDE solutions for both curvatures and all three triangulation types as the resolution is increased. For each resolution, the curves are also almost exactly the same across all three triangulation types, aside from the scalar curvature for the diamond type triangulations.

Since it was shown that the curvature decays exponentially to zero by Carfora et al. \cite{CarfIsenGowdy}, exponential functions were fitted to the time evolution of the scalar and Ricci curvatures graphed in figure \ref{fig:Gowdyt}, with the decay rates for the best-fit functions shown in table \ref{tab:GowdyParam}. Each column can be seen to converge to the PDE value, with good approximations for even the lowest resolutions, and very similar values for each resolution across all three triangulation types, particularly for the Ricci curvature. The $R$-squared values of the fitted functions are greater than $0.99999$ for all but the scalar curvatures of the diamond triangulations, which are no lower than $0.9997$. The larger errors and lower $R$-squared values may be due to the overlapping of edges along the $\theta$-direction for the diamond blocks, giving a larger $\theta$-interval for the averaging of the scalar curvature and requiring more blocks to resolve it better.

\

\begin{table}[h!]
	\centering
	\begin{tabular}{rcc|cc|ccc
		}
		&\multicolumn{2}{c|}{Cubic}
		&\multicolumn{2}{|c|}{Skew}
		&\multicolumn{2}{|c}{Diamond}
		\\
		\cline{2-7}
		& $R$ & $Rc$
		& $R$ & $Rc$
		& $R$ & $Rc$
		\\
		6-vertices \
		& $1.683$ & $0.917$ 
		& $1.732$ & $0.918$
		& $2.874$ & $0.910$
		\\
		12-vertices \
		& $1.869$ & $0.983$ 
		& $1.908$ & $0.984$
		& $1.852$ & $0.981$
		\\
		24-vertices \
		& $1.955$ & $1.001$ 
		& $1.981$ & $1.001$
		& $1.884$ & $1.000$
		\\
		PDE solution \
		& $2.003$ & $1.003$ 
		& $2.003$ & $1.003$
		& $2.003$ & $1.003$
		\\
	\end{tabular}
	\caption{Decay rates for the best-fit exponential functions of the scalar curvature and the Ricci curvature along the $y$-direction, both at $\theta = \pi/3$. The $R$-squared values are at least $0.99999$ for all but the diamond scalar curvatures, and each column shows a convergence to the corresponding PDE value as the resolution is increased.}
	\label{tab:GowdyParam}
\end{table}

\begin{figure}[h!]
	\begin{center}
		\includegraphics[scale=1]{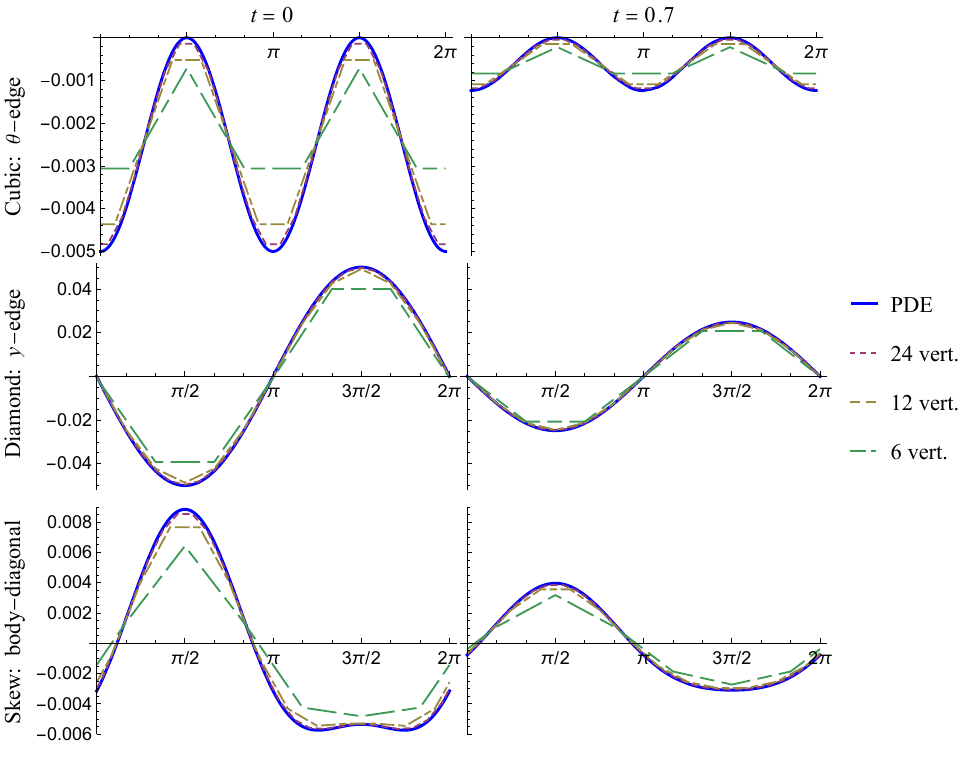}
	\end{center}
	\vspace{-0.5cm}
	\caption{Graphs of the Ricci curvature as a function of $\theta$, for $t=0$ and $t = 0.7$, showing convergence to the PDE solutions as the resolution is increased in all cases.}
	\label{fig:GowdyRc}
\end{figure}

Visual representations of the Ricci curvature as a function of $\theta$ are shown in figure \ref{fig:GowdyRc} for times $t = 0$ and $t = 0.7$, with different edges used for each triangulation type.
The $\theta$-coordinate was selected for the cubic triangulations, since $Rc(\hat \theta, \hat \theta) \equiv R$, corresponding with the scalar curvature in figure \ref{fig:Gowdyt} and table \ref{tab:GowdyParam} above. The $y$-edge was chosen for the diamond triangulation type, corresponding to the Ricci curvature in figure \ref{fig:Gowdyt} and table \ref{tab:GowdyParam}, and the body-diagonal for the skew triangulations.
All six graphs show good approximations for all resolutions, and a convergence to the PDE curves as the resolutions are increased.
The shape of the first two sets of graphs remain mostly unchanged, with the scale decreasing by about a quarter and a half respectively, agreeing with the decay rates in table \ref{tab:GowdyParam}.
Despite the lengths of the body-diagonals being re-defined to give zero deficit angles, required for the stability of the flow as mentioned at the end of section \ref{sec:Tri}, the piecewise flat curves are no further from the PDE curve than for the other edges. This reinforces the robustness of the curvature constructions in (\ref{eq:RandK}). These curves also show the piecewise flat Ricci flow producing more than just a scale change, with the change in shape between the two times resulting from a change in the orientation of the edges, which is perfectly matched by the piecewise flat curvatures.

To show that the edges graphed in figure \ref{fig:GowdyRc} are not special cases, the errors for the Ricci curvatures were computed at all edges for both $t = 0$ and $t = 0.7$, with respect to the corresponding PDE solution values. There is evidence that the errors scale with the overall curvature of the manifold, so the mean absolute values of the errors are given as percentages of an average of the Ricci curvature for the corresponding value of $t$. This average is found by first taking the square-root of the tensor square of the Ricci curvature and then averaging over $\theta$,
\begin{equation}
\label{eq:GRcAve}
Rc_{\textrm ave} := \frac{1}{2 \pi} \int_0^{2 \pi} \sqrt{Rc^2} \, d \theta \, .
\end{equation}
In table \ref{tab:GowdyPercErr}, these percentage errors can be seen to decrease in all cases as the resolution is increased. The errors also have similar values across all triangulation types, showing a certain level of independence to this choice for both the curvature values and the Ricci flow.

\

\begin{table}[h!]
	\centering
	\begin{tabular}{rcc|cc|ccc
		}
		 &\multicolumn{2}{c|}{Cubic}
		 &\multicolumn{2}{|c|}{Skew}
		 &\multicolumn{2}{|c}{Diamond}
		\\
		\cline{2-7}
		& $t = 0$ & $t = 0.7$
		& $t = 0$ & $t = 0.7$
		& $t = 0$ & $t = 0.7$
		\\
		6-vertices \
		& $12.4 \%$ & $6.5 \%$ 
		& $ 9.5 \%$ & $4.5 \%$
		& $ 8.7 \%$ & $3.2 \%$
		\\
		12-vertices \
		&  $3.3 \%$ & $1.3 \%$ 
		&  $2.9 \%$ & $1.1 \%$
		&  $2.3 \%$ & $0.7 \%$
		\\
		24-vertices \
		&  $0.8 \%$ & $0.4 \%$ 
		&  $1.2 \%$ & $1.0 \%$
		&  $0.6 \%$ & $0.4 \%$
		\\
	\end{tabular}
	\caption{The mean absolute error of the Ricci curvature over all edges in each triangulation, as a percentage of the mean absolute value of the curvature at each time (\ref{eq:GRcAve}), showing a decrease everywhere as the resolution is increased.}
	\label{tab:GowdyPercErr}
\end{table}

There is a notable reduction in the errors from $t=0$ to $t=0.7$ for the lower resolution triangulations, the top row of table \ref{tab:GowdyPercErr}. This is likely due to the under-approximation of the initial curvature magnitudes, as seen in figure \ref{fig:GowdyRc}, giving a slower decay rate and reducing these errors in time. The effect indicates a general stabilizing behaviour in the piecewise flat Ricci flow, with errors \emph{reducing} over time, at least for manifolds where the curvature magnitude is decreasing everywhere.

\section{Three-torus initially embedded in $\mathbb{E}^4$}
\label{sec:Torus}

\subsection{Smooth manifold}

\begin{figure}[h!]
	\begin{center}
		\includegraphics[scale=1.2]{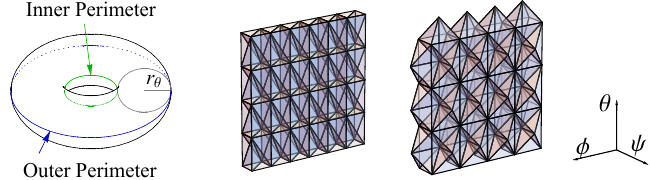}
	\end{center}
	\caption{A two-torus embedded in $\mathbb{E}^3$, showing the inner and outer perimeters tangent to the $\phi$-coordinate vector fields, and the lowest resolution grids for both the cubic and diamond triangulation types, with opposite sides identified to form a $T^3$ topology.}
	\label{fig:T3Tri}
\end{figure}

A regular three-torus can be embedded in Euclidean four-space, $\mathbb{E}^4$, by beginning with a circle of radius $r_\theta$ and rotating it in an orthogonal plane, around a point a distance $r_\phi > r_\theta$ from its centre. The resulting two-torus can then be rotated in a plane orthogonal to this, around a point a distance $r_\psi > r_\phi$ from the centre of the two-torus. For a poloidal-type coordinate system, with angular coordinates $\theta$, $\phi$ and $\psi$ for each rotation, the metric induced by this embedding is
\begin{equation}
\label{eq:T3metric}
d s^2
 = r_\theta^2 \ d \theta^2
 + (r_\phi + \cos \theta)^2 \ d \phi^2
 + (r_\psi +(r_\phi + \cos \theta) \cos \phi)^2 \ d \psi^2 \ .
\end{equation}
This manifold is isometric along the $\psi$-coordinate, and the integral curves of the $\theta$-vector fields are the same everywhere, with values of $2 \pi r_\theta$. The minimum and maximum $\phi$ integral curves have lengths $2 \pi (r_\phi - r_\theta)$ and $2 \pi (r_\phi + r_\theta)$ with $\psi$ integral curves of length $2 \pi (r_\psi - r_\phi)$ to $2 \pi (r_\psi + r_\phi)$. These can be seen as generalizations of the inner and outer perimeters of a two-torus in $\mathbb{E}^3$, as shown in figure \ref{fig:T3Tri}. The smooth Ricci flow of the metric above is expected to tend asymptotically to a flat three-torus, with the minimum and maximum orbits both approaching the same value asymptotically, for each coordinate.

\subsection{Triangulating manifold}

For the piecewise flat computations, the radii were chosen to have the values
\begin{equation}
\label{eq:T3radii}
r_\theta = 1 \, , \qquad r_\phi = 2 \, , \qquad r_\psi = 4 \, .
\end{equation}
Since the manifold is isometric in the $\psi$-direction, only a single layer of blocks in a $\theta \phi$-plane is required. For the cubic block, grids of size $4 \times 6$, $6 \times 6$ and $6 \times 8$ are used, with grids of size $3 \times 4$, $4 \times 4$ and $4 \times 5$ for the diamond block. These give the same number of vertices and edges for the smallest and largest grids for both triangulation types. The entire range of the $\theta$ and $\phi$ coordinates is covered by all triangulations, but to keep the blocks regular in shape, only $1/(2 n_x)$ of the $\psi$-coordinate is covered, where $n_x$ is the number of blocks along the $\theta$-direction. This means that $8$, $12$ and $12$ copies of the cubic triangulations are needed to cover the entire manifold, and $6$, $8$ and $8$ copies of the diamond triangulations. Unfortunately, the skew triangulation would require at least three layers of blocks in the $\psi$-direction to have vertices and edges on opposite sides match, so it has not been used for this manifold.

\subsection{Results of evolution}

\begin{figure}[h!]
	\begin{center}
		\includegraphics[scale=1]{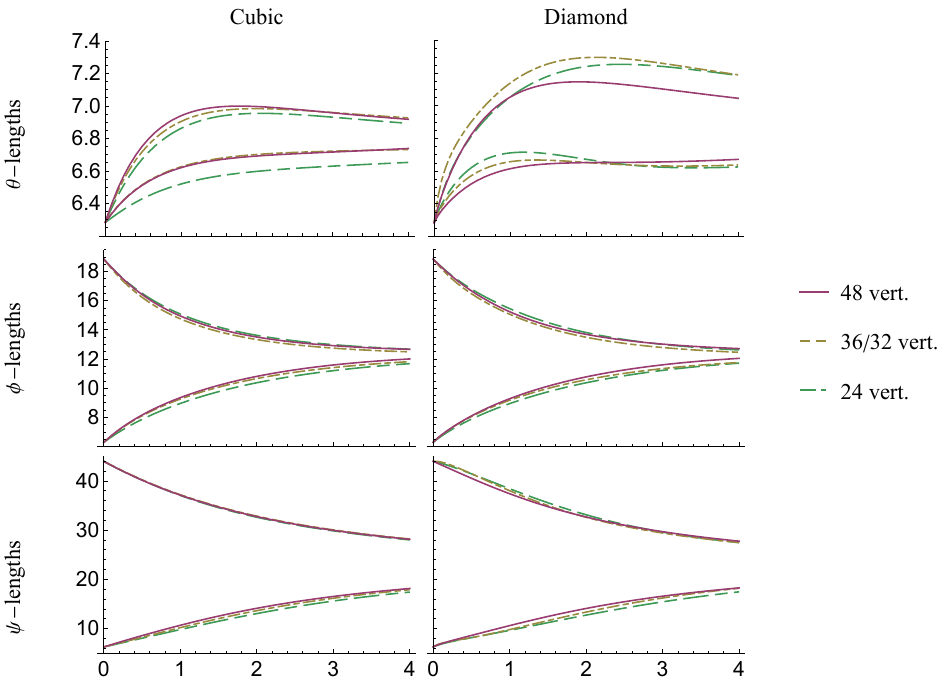}
	\end{center}
	\vspace{-0.5cm}
	\caption{The flow of the maximum and minimum integral curve lengths for each coordinate, tending asymptotically toward the same value in all cases. The curves for both triangulation types also tend to converge to the same shape and values as the resolutions are increased.}
	\label{fig:Torust}
\end{figure}

The triangulations are evolved according to the non-normalized piecewise flat Ricci flow equations (\ref{eq:PFRF}), using the Euler method with $80$ steps of size $0.05$. In order to help visualize the flow, the lengths of the edges along both the minimum and maximum $\phi$ and $\psi$ integral curves are summed for each time step. Although the $\theta$ integral curves all have equivalent lengths to begin with, the Ricci curvatures are not invariant along the $\phi$-coordinate, so different length curves emerge with the flow. The total lengths of the minimum and maximum integral curves in time, for all three resolutions, are graphed together in figure \ref{fig:Torust} for each coordinate and each triangulation type.

\begin{figure}[h!]
	\begin{center}
		\includegraphics[scale=1]{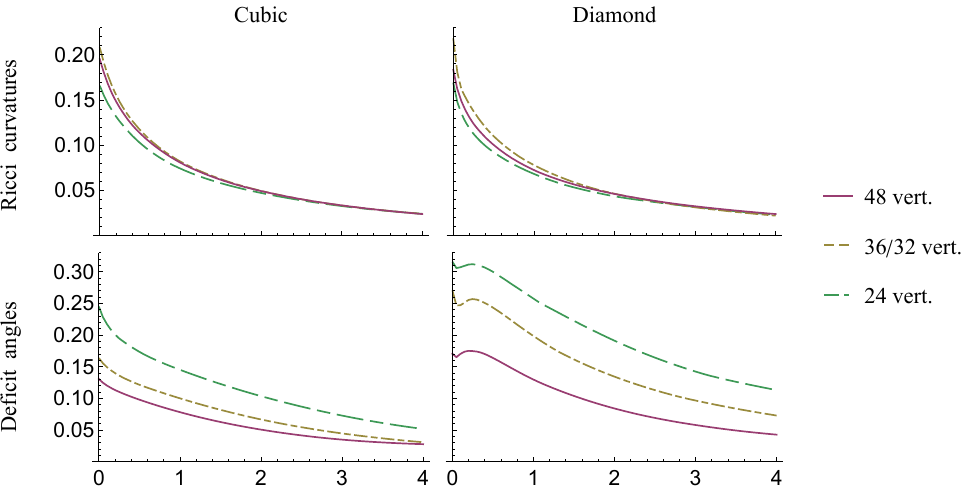}
	\end{center}
	\vspace{-0.5cm}
	\caption{
	The average magnitudes of both the Ricci curvature and deficit angles, all tending asymptotically to zero. The Ricci curvature also tends toward the same curve for both triangulation types, and the deficit angles are smaller for increased resolutions at all times.
	}
	\label{fig:TorusAveRc}
\end{figure}

The graphs show the integral curve lengths converging asymptotically toward the same value for each coordinate, indicating that the manifold is flowing toward a three-torus with consistent dimensions. The values at $t=4$ are also given in table \ref{tab:TorusEnds}, providing bounds for the final three-torus dimensions for each triangulation. Both the curves and table values can also be seen to converge to the same shape and values for both triangulation types as the resolution is increased, sometimes from different directions. This suggests that the higher resolution values are closer to the smooth values, with both triangulation types agreeing to appropriate levels of precision. To show that the manifold is actually flowing to a \emph{flat} three-torus, the average of the absolute values of the Ricci curvature and deficit angles are also computed for each time step. The results are graphed in figure \ref{fig:TorusAveRc}, showing that both flow asymptotically to zero.

\begin{table}[h!]
	\centering
	\begin{tabular}{lcc|cc|cc
		}
		&\multicolumn{2}{c|}{$\theta$}
		&\multicolumn{2}{|c|}{$\phi$}
		&\multicolumn{2}{|c}{$\psi$}
		\\
		\hline
		Cubic \
		& Min & Max
		& Min & Max
		& Min & Max
		\\
		24-Block \ 
		& $6.653$ & $6.893$ 
		& $11.68$ & $12.67$ 
		& $17.43$ & $28.02$ 
		\\
		36-Block \ 
		& $6.732$ & $6.927$ 
		& $11.82$ & $12.49$ 
		& $17.91$ & $28.21$ 
		\\
		48-Block \ 
		& $6.738$ & $6.918$ 
		& $12.01$ & $12.66$ 
		& $18.16$ & $28.15$ 
		\\
		&&&&& \\
		\hline
		Diamond \
		& Min & Max
		& Min & Max
		& Min & Max
		\\
		12-Block \ 
		& $6.626$ & $7.187$ 
		& $11.72$ & $12.62$ 
		& $17.48$ & $27.47$ 
		\\
		16-Block \ 
		& $6.638$ & $7.190$ 
		& $11.76$ & $12.46$ 
		& $18.21$ & $27.48$ 
		\\
		24-Block \ 
		& $6.672$ & $7.046$ 
		& $12.04$ & $12.72$ 
		& $18.28$ & $27.76$ 
		\\
	\end{tabular}
	\caption{
	The minimum and maximum perimeter lengths at $t=4$, bounding the dimensions of the limiting flat three-torus, with both triangulation types tending towards the same values.	
	}
	\label{tab:TorusEnds}
\end{table}

\section{Perturbation of a flat three-torus}
\label{sec:Pert}

\subsection{Smooth manifold}

A simple perturbation of a flat three-torus has been used to give a manifold without any continuous isometries. For a three-torus topology, with coordinates $x$, $y$ and $z$ ranging from $0$ to $1$, the metric was chosen to be
\begin{equation}
\label{eq:PertMetric}
d s^2 = (1 + 0.2 \sin^2 (\pi x) \sin^2 (\pi y) \sin^2 (\pi z))\left(d x^2 + d y^2 + d z^2\right) \ .
\end{equation}
As with both the Gowdy manifold and the three-torus, this manifold is expected to Ricci flow asymptotically back to a flat three-torus.

\subsection{Triangulating manifold}

Without any continuous isometries, there are no reductions that can be made for the triangulation grids. The cubic triangulations have blocks arranged in rectangular grids with $2$, $3$ and $4$ blocks in each direction, while the diamond triangulations start with a single block, with grids of $2$ and $3$ blocks in each direction as well. The vertices and edges of the skew blocks must align with other vertices and edges when opposite faces are identified. To make this easier for smaller grid sizes, the skew block is adapted a little with the vertices in the $x$ and $z$ directions skewed so that $v_x = (1, -1/2, 0)$ and $v_z = (-1/2, -1/4, 1)$. This still requires even numbers of blocks in each direction, so only grids of $2$ and $4$ blocks in each direction are used.

\begin{figure}[h!]
	\begin{center}
		\includegraphics[scale=1.2]{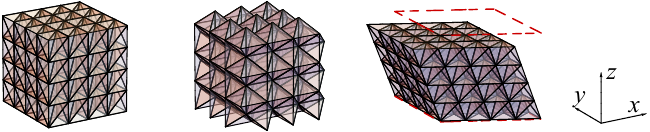}
	\end{center}
	\caption{
	The highest resolution grids for the cubic, diamond and skew triangulation types, with the identification between the two $x y$-faces indicated for the skew triangulation.
	}
	\label{fig:PertTri}
\end{figure}

\subsection{Results of evolution}

\begin{figure}[h!]
	\begin{center}
		\includegraphics[scale=1]{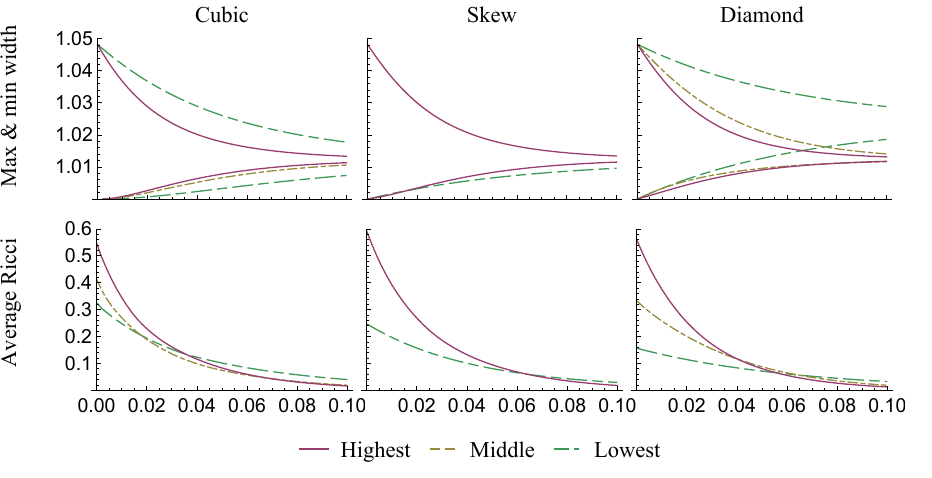}
	\end{center}
	\caption{
	On the top row, the minimum and maximum integral curve lengths along the $y$-coordinate show convergence to the same shape for each triangulation type. The graphs on the bottom row show the curvature flowing asymptotically to zero.
	}
	\label{fig:Pertt}
\end{figure}

The triangulations are evolved using the Euler method, with $50$ steps of size $0.002$, for the non-normalized piecewise flat Ricci flow. The smooth manifold is expected to flow asymptotically to a flat three-torus, so the lengths of the integral curves for each coordinate should all converge to the same value. Though the manifold does not have any continuous isometries, it is discretely symmetric for any permutation of the coordinates, so the three sets of coordinate integral curves should be equivalent. 
Since all of the triangulation types have edges aligned with the $y$-coordinate, the minimum length occurring at $(x, z) = (0, 0)$ and maximum length at $(x, z) = (0.5, 0.5)$ are computed and graphed in figure \ref{fig:Pertt}. While the minimum lengths can be computed for all of the triangulations, there are no $y$-edges at $(x, z) = (0.5, 0.5)$ for the middle resolution cubic triangulation and lowest resolution skew triangulation, so these are omitted from the graph. The average magnitude of the Ricci curvatures at each time-step are also graphed for each triangulation on the bottom row of figure \ref{fig:Pertt}.

The bottom row of figure \ref{fig:Pertt} shows the Ricci curvature asymptotically approaching zero for all of the triangulations, with a consistent shape to the curve for the highest resolutions. This indicates that the manifold is flowing towards a flat three-torus. 
The curves for the minimum and maximum lengths can also be seen to tend asymptotically to the same value for each of the triangulation types, indicating that the manifold tends to a three-torus with consistent dimensions. The values of these lengths at $t = 0.1$ are also given in table \ref{tab:PertEnds}, giving bounds for the dimensions of the limiting manifold, with both the minimum and maximum values converging towards similar values for all three triangulation types as the resolution is increased.

\begin{table}[h!]
	\centering
	\begin{tabular}{lcc|cc|cc
		}
		&\multicolumn{2}{c|}{Cubic}
		&\multicolumn{2}{|c|}{Skew}
		&\multicolumn{2}{|c}{Diamond}
		\\
		\hline
		Resolutions
		& Min & Max
		& Min & Max
		& Min & Max
		\\
		Lowest 
		& $1.0074$ & $1.0178$ 
		& $1.0095$ & $-$ 
		& $1.0186$ & $1.0288$ 
		\\
		Middle 
		& $1.0106$ & $-$ 
		& $-$ & $-$ 
		& $1.0117$ & $1.0141$ 
		\\
		Highest 
		& $1.0113$ & $1.0134$ 
		& $1.0115$ & $1.0134$ 
		& $1.0118$ & $1.0132$ 
		\\
	\end{tabular}
	\caption{
	Values for both the minimum and maximum lengths in the $y$-direction at $t=0.1$, giving bounds on the dimensions of the limiting flat three-torus for each approximation.
	}
	\label{tab:PertEnds}
\end{table}

\section{Conclusion}
\label{sec:Con}

The computations in sections \ref{sec:Nil} to \ref{sec:Pert} have successfully demonstrated the effectiveness of the piecewise flat Ricci flow introduced in \cite{PLCurv}, showing
\begin{itemize}
	\item a clear convergence to the smooth Ricci flow as the mesh is refined, 
	
	\item giving equivalent results for different mesh types, up to appropriate levels of precision,
	
	\item for manifolds varying from completely homogeneous to completely inhomogeneous.
	
\end{itemize}

The convergence to the expected smooth Ricci flow has been shown in a number of ways. In the case of the Nil and Gowdy manifolds, where smooth Ricci flow solutions exist, appropriate equations for the metric components and curvature values give extremely close fits to the evolution data, with $R$-squared values that differ from 1 only at the sixth or seventh decimals. Tables \ref{tab:NilParam} and \ref{tab:GowdyParam} then clearly show the parameters for these fits converging to their corresponding smooth Ricci flow values as the mesh resolution increases. The remaining two manifolds are each expected to Ricci flow asymptotically to a flat three torus. For the piecewise flat Ricci flow, figures \ref{fig:TorusAveRc} and \ref{fig:Pertt} show the average magnitude of the Ricci curvatures approaching zero asymptotically, while the dimensions of the limiting domain are shown to be consistent with minimum and maximum widths in each direction approaching the same values, see figures \ref{fig:Torust} and \ref{fig:Pertt}, and tables \ref{tab:TorusEnds} and \ref{tab:PertEnds}.

For each manifold, the different triangulations also converge to consistent evolutions. 
This is particularly clear for edge orientations that are common across different triangulation types. For the Gowdy manifold, the best-fit parameters for a $y$-edge all approach the same smooth solution values in table \ref{tab:GowdyParam}. The evolutions for the minimum and maximum widths in figures \ref{fig:Torust} and \ref{fig:Pertt} converge toward consistent shaped curves for all triangulation types and, interestingly, the graphs for the $\theta$-lengths on the top row of the former converge from opposite directions. This implies that different triangulations may be able to provide upper and lower bounds for the limiting smooth Ricci flow manifold. 
The computations also show that the triangulations do not need to be adapted to any manifold symmetries for the evolutions to work successfully, as other approaches do. While the cubic and diamond triangulation types do have edges that align with all three coordinate base vectors, the skew triangulation does not, and the Nil metric is not even diagonal for these base vectors.

\section*{Declarations}

\noindent
\emph{Data Availability:} The \emph{Mathematica} notebooks used for the computations and numerical simulations, and the data generated by these, are available in the Zenodo repository at https://doi.org/10.5281/zenodo.10697618.

\

\noindent
\emph{Conflict of Interest:} The author has no competing interests to declare.

\

\noindent
\emph{Funding:} No funds, grants, or other support was received.


\appendix

\bibliography{SimpGeom}
\bibliographystyle{unsrt}

\end{document}